\documentclass[journal]{IEEEtran}

\usepackage[noadjust]{cite}
\usepackage{amsmath,amssymb,amsfonts,mathtools}
\usepackage{algorithmic}
\usepackage{arydshln}
\usepackage{graphicx}

\usepackage[inline,shortlabels]{enumitem}

\usepackage{amsthm}
\usepackage{eso-pic}
\AddToShipoutPictureBG*{%
\AtPageUpperLeft{%
\setlength\unitlength{1in}%
\hspace*{\dimexpr0.5\paperwidth\relax}
\makebox(0,-1)[c]{
\begin{tabular}{c c}
T.~J.~{Meijer} \emph{et al.},
The Non-Strict Projection Lemma, \\
To appear in {\em IEEE Transactions on Automatic Control}, uploaded to ArXiv \today \\
\end{tabular}}}}
\AddToShipoutPictureBG*{%
\AtPageUpperLeft{%
\setlength\unitlength{1in}%
\hspace*{\dimexpr0.5\paperwidth\relax}
\makebox(0,-21.4)[c]{
\footnotesize
\begin{tabular}{c c}
© 2024 IEEE. Personal use of this material is permitted. Permission from
IEEE must be obtained for all other uses, in any \\
current or future media, including reprinting/republishing
this material for advertising or promotional purposes, creating new \\
collective works, for resale or redistribution to servers or lists,
or reuse of any copyrighted component of this work in other works.
\end{tabular}}}}

\newcommand*{\suchthat}{%
    \;%
    \ifnum\currentgrouptype=16\relax
        \middle|%
    \else%
        \iftoggle{WithinBracMacro}{
            |
        }{
            |%
        }%
    \fi%
    \;%
}%
\makeatother
\newcommand{\im}{\operatorname{im} }
\newcommand{\hop}{\mathsf{H} }

\newtheorem{lemma}{Lemma}
\newtheorem{example}{Example}
\newtheorem{theorem}{Theorem}
\newtheorem{proposition}{Proposition}
\newtheorem{corollary}{Corollary}
\newtheorem{remark}{Remark}
\theoremstyle{definition}

\makeatletter
\newcommand\footnoteref[1]{\protected@xdef\@thefnmark{\ref{#1}}\@footnotemark}
\makeatother

\DeclareRobustCommand{\svdots}{
  \vbox{%
    \baselineskip=0.33333\normalbaselineskip
    \lineskiplimit=0pt
    \hbox{.}\hbox{.}\hbox{.}%
    \kern-0.2\baselineskip
  }%
}

\begin{document}
\title{The Non-Strict Projection Lemma}
\author{T.~J.~{Meijer}, T.~{Holicki}, S.~J.~A.~M.~{van den Eijnden},\\ C.~W.~{Scherer},~\IEEEmembership{Fellow, IEEE}, and W.~P.~M.~H.~{Heemels},~\IEEEmembership{Fellow, IEEE}
\thanks{{Tomas} {Meijer}, {Sebastiaan} {van den Eijnden} and {Maurice Heemels} are with the Department of Mechanical Engineering, Eindhoven University of Technology, the Netherlands. (e-mail: t.j.meijer@tue.nl; s.j.a.m.v.d.eijnden@tue.nl; m.heemels@tue.nl).}
\thanks{{Tobias} {Holicki} and {Carsten} {Scherer} are with the Department of Mathematics, University of Stuttgart, Germany (e-mail: tobias.holicki@mathematik.uni-stuttgart.de; carsten.scherer@imng.uni-stuttgart.de)}%
\thanks{This research received funding from the European Research Council
(ERC) under the Advanced ERC grant agreement PROACTHIS, no.
101055384.}%
\thanks{The second and fourth author are funded by Deutsche Forschungsgemeinschaft (DFG, German Research Foundation) under Germany’s Excellence Strategy -EXC 2075 -390740016. They acknowledge the support by the Stuttgart Center for Simulation Science (SimTech).}%
\thanks{Corresponding author: T.~J.~{Meijer}.}}

\maketitle

\begin{abstract}
	The projection lemma (often also referred to as the elimination lemma) is one of the most powerful and useful tools in the context of linear matrix inequalities for system analysis and control. In its traditional formulation, the projection lemma only applies to strict inequalities, however, in many applications we naturally encounter non-strict inequalities. As such, we present, in this note, a non-strict projection lemma that generalizes both its original strict formulation as well as an earlier non-strict version. We demonstrate several applications of our result in robust linear-matrix-inequality-based marginal stability analysis and stabilization, a matrix S-lemma, which is useful in (direct) data-driven control applications, and matrix dilation.
\end{abstract}


\section{Introduction}
\label{sec:introduction}
\IEEEPARstart{L}{inear} matrix inequalities (LMIs) have found their way into a wide variety of control applications \cite{Boyd1994}. In parallel to this adoption, an incredible collection of tools has been developed that enables us to formulate LMIs for different and increasingly complicated applications. The projection lemma (PL), see, e.g.,~\cite{Gahinet1994,Scherer2000}, is a crucial part of this LMI toolkit, which has been an enabler for developing powerful results, such as, e.g., $H_\infty$ controller synthesis~\cite{Gahinet1994,Feng2013}, robust control design~\cite{Packard1991, Boyd1994, Skelton1998}, and gain-scheduled control design \cite{Packard1994, Apkarian1995} to name but a few. In fact,~\cite[Chapter 9]{Skelton1998} presents a unified approach based on the projection lemma to solve 17 seemingly different control problems, including the characterization of all stabilizing controllers for a linear time-invariant (LTI) plant, covariance control, $H_\infty$ control, $L_\infty$ control, linear-quadratic-Gaussian (LQG) control, and $H_2$ control of LTI systems. It is also used to introduce slack variables for reducing the conservatism in certain robust control designs, see, e.g.,~\cite{Oliveira2002,Ebihara2015}. In other words, the PL has had--without a doubt--a significant impact in the field of system and control theory. All of the above developments are based on the strict version of the projection lemma (strict in the sense of strictness of the involved matrix inequalities). Given the impact of this strict projection lemma (SPL) and the emergence of various control problems that call for non-strict versions of the PL (see Sections~\ref{sec:motivating-example} and~\ref{sec:applications} below), we will formulate a non-strict generalization of this powerful result.

The classical (strict) projection lemma~\cite{Boyd1994,Gahinet1994}, as stated next, is formulated in terms of strict inequalities. For any complex matrix $A\in\mathbb{C}^{m\times n}$, we denote its conjugate transpose by
$A^\hop$ and its annihilator by $A_\perp$, which is any matrix whose columns form a basis of the kernel (null space) of $A$.

\begin{lemma}\label{lem:SPL}
    Let $U\in\mathbb{C}^{m\times p}$ and $V\in\mathbb{C}^{n\times p}$
    be arbitrary complex matrices and let $Q\in\mathbb{H}^{p}$ be Hermitian.
    Then, there exists a matrix $X\in\mathbb{C}^{m\times n}$ which satisfies the LMI
    \begin{equation}
        Q + U^\hop XV + V^\hop X^\hop U\succ 0,
        \label{eq:trad-proj-1}
    \end{equation}
    if and only if
    \begin{equation}
        U^\hop_{\perp}Q U_{\perp}\succ 0\text{ and }V^\hop_{\perp}Q V_{\perp}\succ 0.
        \label{eq:trad-proj-2}
    \end{equation}
\end{lemma}

\noindent
When replacing the \emph{strict} inequalities ($\succ$) in Lemma~\ref{lem:SPL} by \emph{non-strict} inequalities ($\succcurlyeq$), the implication \eqref{eq:trad-proj-1} $\Rightarrow$ \eqref{eq:trad-proj-2} still holds.
However, the converse is no longer true as illustrated with the following example.
\begin{example}\label{ex-non-sufficiency}
    Consider the matrices
    \begin{equation*}
        Q = \begin{bmatrix}
            2 & 1 \\
            1 & 0
        \end{bmatrix},~U = V = \begin{bmatrix} 1 & 0\end{bmatrix}\text{ with }U_{\perp}=V_{\perp}=\begin{bmatrix}0\\ 1\end{bmatrix}
    \end{equation*}
    as the relevant annihilators. It is straightforward to verify that $U_\perp^\hop QU_\perp=V^\hop_\perp QV_\perp = 0\succcurlyeq 0$, i.e., the non-strict version of~\eqref{eq:trad-proj-2} holds. However, there  exists no $X\in\mathbb{C}$ such that
    \begin{equation*}
        Q+U^\hop XV + V^\hop X^\hop U = \begin{bmatrix} 2 + (X + X^*) & 1\\ 1 & 0\end{bmatrix}\succcurlyeq 0,
    \end{equation*}
    where $X^*$ denotes the complex conjugate of $X$. This shows that the non-strict version of~\eqref{eq:trad-proj-2} alone is not sufficient to guarantee the existence of some $X$ for which~\eqref{eq:trad-proj-1} holds with a non-strict inequality.
\end{example}
\noindent The following non-strict version of the PL
from~\cite[Lemma 6.3]{Helmersson1995} imposes a strong additional assumption on $U$ and $V$ under which~\eqref{eq:trad-proj-1}~$\Leftrightarrow$~\eqref{eq:trad-proj-2} holds with non-strict inequalities as well.
\begin{lemma}\label{lem:ns-proj-Helmersson}
	Suppose that $U$ and $V$ satisfy
	\begin{equation}
		\im U^\hop \cap\im V^\hop =\{0\}.
		\label{eq:Helmersson-extra-cond}
	\end{equation}
	Then, there exists some $X$ such that
	\begin{equation}\label{eq:ns-proj-lem-X}
		Q + U^\hop XV+V^\hop X^\hop U\succcurlyeq 0,
	\end{equation}
	if and only if
	\begin{equation}\label{eq:ns-proj-lem-kernels}
		U_{\perp}^\hop Q U_{\perp}\succcurlyeq 0\text{ and }V_{\perp}^{\hop }Q V_{\perp}\succcurlyeq 0.
	\end{equation}
\end{lemma}
\noindent
This result, in which $\operatorname{im}\left(\cdot\right)$ denotes the image of a matrix,
is shown to be useful in the context of, e.g., robust control using $\mu$-synthesis~\cite{Helmersson1995}.
However, due to the strong assumption \eqref{eq:Helmersson-extra-cond}, this result cannot be employed in many situations, such as for characterizing marginal Lyapunov stability as we will illustrate in Section~\ref{sec:motivating-example}. Moreover, the condition~\eqref{eq:Helmersson-extra-cond} also prevents Lemma~\ref{lem:ns-proj-Helmersson} from yielding the SPL as a corollary, since the SPL does not involve any assumptions on the matrices $U$ and $V$. In this sense, Lemma 2 is not a true generalization of the SPL.

Our main contribution in this note is a general non-strict projection lemma (NSPL), see Theorem~\ref{thm:ns-proj-lem} below, that
\begin{enumerate}[label=(\alph*)]
	\item provides necessary and sufficient conditions for the existence of a solution to the LMI \eqref{eq:ns-proj-lem-X},
	
	\item does not impose \emph{any} assumption on the matrices $U$, $V$ and $Q=Q^\hop$, and
	
	\item generalizes both the SPL and Lemma~\ref{lem:ns-proj-Helmersson}.
\end{enumerate}
	
\noindent Although special cases of the NSPL have found applications in earlier work, see, e.g.,~\cite{Helmersson1995,Scherer2005,Scherer2000}, the general NSPL is, to the best of the authors' knowledge, not found in the literature.

One of the key challenges in deriving a general NSPL is the identification of an additional condition that, if combined with~\eqref{eq:ns-proj-lem-kernels}, is equivalent to the feasibility of~\eqref{eq:ns-proj-lem-X}. Moreover, even though our proof of the NSPL is inspired by the proof of the SPL, the non-strictness of the inequalities and the absence of simplifying assumptions, such as~\eqref{eq:Helmersson-extra-cond}, result in additional technical challenges. These challenges are highlighted in detail in our proof of the NSPL.

We illustrate how to apply our main result for several applications. First, we derive LMI-based marginal stability conditions, which could not be achieved using the SPL or Lemma~\ref{lem:ns-proj-Helmersson}.
Second, we apply the NSPL to solve a matrix dilation problem and, third, we derive a useful interpolation result with weaker assumptions than existing ones. Interestingly, this interpolation result naturally leads to a generalization of the matrix S-lemma in~\cite[Corollary 12]{vanWaarde2022}, which is used for (direct) data-driven control.

The remainder of this note is organized as follows. After introducing some notational conventions in Section~\ref{sec:notation}, we discuss a motivating example, in Section~\ref{sec:motivating-example}, for which the SPL and Lemma~\ref{lem:ns-proj-Helmersson} fall short. In Section~\ref{sec:main-result}, we present the NSPL, which forms our main contribution. Finally, we present several applications of this NSPL in Section~\ref{sec:applications} and give conclusions in Section~\ref{sec:conclusions}. All proofs are found in the Appendix.

\subsection{Notation}\label{sec:notation}
The sets of real, complex and non-negative natural numbers are denoted, respectively, by $\mathbb{R}$, $\mathbb{C}$ and $\mathbb{N}=\{0,1,2,\hdots\}$. The sets of $n$-dimensional real and complex vectors are, respectively, $\mathbb{R}^{n}$ and $\mathbb{C}^{n}$. The sets of $n$-by-$n$ Hermitian and symmetric matrices are denoted by, respectively, $\mathbb{H}^{n}=\{A\in\mathbb{C}^{n\times n}\mid A=A^\hop\}$ and $\mathbb{S}^{n}=\{A\in\mathbb{R}^{n\times n}\mid A=A^\top\}$. We use the symbol $\star$ to complete a Hermitian matrix, e.g., $\begin{bmatrix}\begin{smallmatrix} A & B\\\star & C\end{smallmatrix}\end{bmatrix}=\begin{bmatrix}\begin{smallmatrix}A & B\\B^{\hop } & C\end{smallmatrix}\end{bmatrix}$, and $I$ is an identity matrix of appropriate dimension. For a Hermitian matrix $H\in\mathbb{H}^n$, $H\succ 0$, $H\succcurlyeq 0$ and $H\prec 0$ mean, respectively, that $H$ is positive definite, i.e., $x^\hop Hx>0$ for all $x\in\mathbb{C}^{n}\setminus\{0\}$, positive semi-definite, i.e., $x^\hop Hx\geqslant 0$ for all $x\in\mathbb{C}^{n}$, and negative definite, i.e., $-H\succ 0$. We denote, respectively, the sets of such matrices of size $n$-by-$n$ as $\mathbb{H}_{\succ 0}^{n}$, $\mathbb{H}_{\succcurlyeq 0}^{n}$ and $\mathbb{H}_{\prec 0}^{n}$, and their real-valued counterparts as $\mathbb{S}_{\succ 0}^{n}$, $\mathbb{S}_{\succcurlyeq 0}^{n}$ and $\mathbb{S}_{\prec 0}^{n}$. For a complex matrix $A\in\mathbb{C}^{n\times m}$, $\im A=\{x\in\mathbb{C}^{n}\mid x=Ay\text{ for some }y\in\mathbb{C}^{m}\}$ denotes its image, $\ker A=\{x\in\mathbb{C}^{m}\mid Ax=0\}$ its kernel and $A^+$ its (Moore-Penrose) pseudoinverse. Finally, $\operatorname{diag}\{A_1,A_2,\hdots,A_n\}$ denotes a block-diagonal matrix with diagonal blocks $A_i$, $i\in\{1,2,\hdots,n\}$.

\section{Motivating example}\label{sec:motivating-example}
Consider the discrete-time LTI system
\begin{equation}
    x_{k+1} = Ax_k,
    \label{eq:system}
\end{equation}
where $x_k\in\mathbb{R}^{n_x}$ denotes the state at time $k\in\mathbb{N}$. We are interested in marginal stability of the system~\eqref{eq:system}, i.e., whether, for all $x_0$, the solution $x_k$ is uniformly bounded in the sense that there exists $c\geqslant 0$ such that $\|x_k\|\leqslant c\|x_0\|$ for all $k\in\mathbb{N}$. The system~\eqref{eq:system} is marginally stable, if and only if there exists a symmetric matrix $P\in\mathbb{S}^{n_x}$ such that~\cite[p.~211]{Saberi2012}
\begin{equation}
    P\succ 0\text{ and }P-A^\top PA\succcurlyeq 0,\label{eq:stab-conds-1}
\end{equation}
or, equivalently, if there exists $S\in\mathbb{S}^{n_x}$ such that
\begin{equation}
    S\succ 0\text{ and }S-ASA^\top \succcurlyeq 0,
    \label{eq:marg-stab-lyap2}
\end{equation}
in which case $V(x) = x^\top Px$ (with $P=S^{-1}$) is a weak Lyapunov function, i.e., $V$ is positive definite, radially unbounded and non-increasing along solutions to~\eqref{eq:system}~\cite{Kellett2004,Khalil1996}.

While~\eqref{eq:marg-stab-lyap2} can be used to guarantee marginal stability of~\eqref{eq:system}, it cannot easily be extended to more complicated applications such as synthesis or even robust control, e.g., when $A$ is uncertain, by applying~\eqref{eq:marg-stab-lyap2} to the relevant closed-loop dynamics. In the context of \textit{asymptotic} stability, many different LMI-based conditions have been proposed that not only guarantee that~\eqref{eq:system} is asymptotically stable, but also accommodate such more complicated applications, see, e.g.,~\cite{Oliveira1999,Heemels2010,Daafouz2001,Oliveira2001}. Many of these results, see, e.g.,~\cite{Oliveira1999,Daafouz2001,Heemels2010} (or~\cite{Apkarian2001} for continuous time), follow by application of the SPL. For illustrative purposes, we consider the following condition: The system~\eqref{eq:system} is asymptotically stable if and only if there exist a symmetric matrix $S\in\mathbb{S}^{n_x}$ and a matrix $X\in\mathbb{R}^{n_x\times n_x}$ such that~\cite{Oliveira2002}
\begin{equation}
    \begin{bmatrix}
        S & AX\\
        \star & X+X^\top - S
    \end{bmatrix}\succ 0.
    \label{eq:example-strict-cond}
\end{equation}
This follows from the conditions $S\succ 0$ and $S-ASA^\top\succ 0$, which are necessary and sufficient for asymptotic stability, by applying Lemma~\ref{lem:SPL}. To see this, note that~\eqref{eq:example-strict-cond} is~\eqref{eq:trad-proj-1} with
\begin{equation*}
    Q = \begin{bmatrix}
        S & 0\\
        \star & -S
    \end{bmatrix},~U^\top = \begin{bmatrix}
        0\\
        I
    \end{bmatrix}\text{ and }V^\top = \begin{bmatrix}
        A\\
        I
    \end{bmatrix}.
\end{equation*}
This condition is useful for stabilizing controller synthesis by replacing $A$ with $A+BK$ and applying the linearizing change of variables $X = KY$, where $K$ is the to-be-designed controller gain. What makes the condition in~\eqref{eq:example-strict-cond} even more powerful is the absence of products between $A$ and $S$, which enables a natural extension to robust controller synthesis as presented in, e.g.,~\cite{Daafouz2001,Oliveira1999}, in which both $A$ and $S$ depend on some uncertain parameter. Clearly, we cannot use the SPL to obtain a non-strict counterpart to~\eqref{eq:example-strict-cond} that is also equivalent to~\eqref{eq:marg-stab-lyap2}. Also the condition in~\eqref{eq:Helmersson-extra-cond} is not always satisfied and, hence, we cannot apply Lemma~\ref{lem:ns-proj-Helmersson} either. To see this, observe that, if there exists some $x\in\mathbb{R}^{n_x}\setminus\{0\}$ with $Ax=0$, then $(0,x^\top)^\top$ is contained in both $\im U^\top$ and $\im V^\top$ and, hence, their intersection is non-trivial. In the next section, we will present a generalization of Lemma~\ref{lem:ns-proj-Helmersson} after which we will revisit this example in Section~\ref{sec:marginal-stability-revisited} and derive a non-strict version of~\eqref{eq:example-strict-cond}.

\section{Main result}\label{sec:main-result}
The main contribution of this note, i.e., a non-strict generalization of the well-known SPL, is stated below.
\begin{theorem}\label{thm:ns-proj-lem}
    There exists $X$ satisfying \eqref{eq:ns-proj-lem-X} if and only if \eqref{eq:ns-proj-lem-kernels} holds together with
    \begin{equation}
        \ker U \cap \ker V \cap \{\xi\in\mathbb{C}^{p}\mid \xi^\hop Q \xi=0\}\subset \ker Q.
        \label{eq:ns-proj-coupling}
    \end{equation}
\end{theorem}
\noindent We emphasize that Theorem~\ref{thm:ns-proj-lem} does, in contrast to Lemma~\ref{lem:ns-proj-Helmersson}, not involve any additional assumption on the matrices $U$ and $V$ and, thereby, genuinely provides necessary and sufficient conditions for the feasibility of \eqref{eq:ns-proj-lem-X}.
Furthermore, the novel condition \eqref{eq:ns-proj-coupling} does, in contrast to \eqref{eq:ns-proj-lem-X} and \eqref{eq:ns-proj-lem-kernels}, not constitute an LMI. Instead, it is a nontrivial coupling condition between the matrices $U$, $V$ and $Q$ that is absent in the SPL. In the remainder of this section and in the next section, we consider several important applications of Theorem 1 and we illustrate how to verify whether~\eqref{eq:ns-proj-coupling} holds in the related proofs.
In particular, we will see that most of these applications do not exhibit \eqref{eq:ns-proj-coupling} explicitly, since this coupling condition will be shown to be satisfied automatically as a consequence of the underlying problem structure.

\begin{example}
		In Example \ref{ex-non-sufficiency} we have seen that, for the given matrices $U$, $V$ and $Q$, the LMI \eqref{eq:ns-proj-lem-X} is not feasible in $X$ even though \eqref{eq:ns-proj-lem-kernels} holds. By Theorem~\ref{thm:ns-proj-lem}, we infer that \eqref{eq:ns-proj-coupling} must fail to hold. Indeed, for this example, we have $\ker Q=\{0\}$ and $\ker U\cap\ker V\cap \{\xi\in\mathbb{C}^{p}\mid\xi^\hop Q\xi=0\}=\im U_\perp\not\subset\ker Q$.
\end{example}

By using a perturbation argument, we can easily prove the SPL in Lemma~\ref{lem:SPL} by using the NSPL in Theorem~\ref{thm:ns-proj-lem}.

\begin{corollary}\label{cor:SPL}
	The SPL is a special case of the NSPL.
\end{corollary}

We can also show Lemma~\ref{lem:ns-proj-Helmersson} using the NSPL in Theorem~\ref{thm:ns-proj-lem}.

\begin{corollary}\label{cor:Helmersson}
	Lemma~\ref{lem:ns-proj-Helmersson} is a special case of Theorem~\ref{thm:ns-proj-lem}.
\end{corollary}
\noindent To illustrate that Theorem~\ref{thm:ns-proj-lem} is indeed more general than Lemma~\ref{lem:ns-proj-Helmersson} and to make the conservatism introduced by~\eqref{eq:Helmersson-extra-cond} in Lemma~\ref{lem:ns-proj-Helmersson} more concrete, we provide an algebraic example.

\begin{example}\label{ex:generality}
Consider the matrices
\begin{equation*}
    Q = \begin{bmatrix}
        3 & 1 & -2\\
        \star & 1 & -1\\
        \star & \star & 1
    \end{bmatrix},~U^\top = \begin{bmatrix}
        1 & 0\\
        1 & 1\\
        0 & 1
    \end{bmatrix}\text{ and }V^\top = \begin{bmatrix}
        1\\
        0\\
        -1
    \end{bmatrix}.
\end{equation*}
Relevant annihilators are given by
\begin{equation*}
    U_\perp = \begin{bmatrix}
        1 & -1 & 1
    \end{bmatrix}^\top \text{and }V_{\perp} = \begin{bmatrix}
        1 & -1 & 1\\
        0 & 1 & 0
    \end{bmatrix}^\top.
\end{equation*}
It is straightforward to verify that
\begin{equation*}
    U_{\perp}^\top Q U_{\perp}=1\succcurlyeq 0\text{ and }V_{\perp}^\top Q V_{\perp}=\begin{bmatrix} 1 & -1\\\star & 1\end{bmatrix}\succcurlyeq 0.
\end{equation*}
Since $\operatorname{im} U^\top \cap \operatorname{im} V^\top = \im\begin{bmatrix}
        1 & 0 & -1
\end{bmatrix}^\top\neq \{0\}$, we infer that~\eqref{eq:Helmersson-extra-cond} does not hold. Hence, we cannot apply Lemma~\ref{lem:ns-proj-Helmersson}. However, we can still use Theorem~\ref{thm:ns-proj-lem} to conclude that there exists a solution $X$ to~\eqref{eq:ns-proj-lem-X}, because~\eqref{eq:ns-proj-coupling} is trivially satisfied. To see this, note that $U_{\perp}^\top QU_{\perp}=1\succ 0$, such that there does not exist a nonzero $x\in\ker U\cap\ker V$ with $x^\top Q x=0$ at all.
\end{example}

\begin{remark}
	The proof of Theorem~\ref{thm:ns-proj-lem} is constructive.
It also shows that, if $U$, $V$ and $Q$ are real-valued, the result holds with the solution $X$ to~\eqref{eq:ns-proj-lem-X} being real-valued as well.
\end{remark}

Although a rigorous proof has not been published before, special versions of the NSPL have already proved useful in formulating LMI relaxations in robust control~\cite{Scherer2005}. In the next section, we revisit the motivating example from Section~\ref{sec:motivating-example}, for which we showed that the SPL and Lemma~\ref{lem:ns-proj-Helmersson} did not apply, and utilize Theorem~\ref{thm:ns-proj-lem} to find a solution. We will also demonstrate applications of Theorem~\ref{thm:ns-proj-lem} to matrix dilation theory, interpolation and a matrix S-lemma, which is used in modern data-driven techniques.

\section{Applications}\label{sec:applications}
In this section, we present several relevant applications of Theorem~\ref{thm:ns-proj-lem}. For each of the applications below, the crucial steps in their respective proofs are achieved through the NSPL.

\subsection{Marginal stability and stabilizability revisited}\label{sec:marginal-stability-revisited}
First, we revisit the motivating example discussed in Section~\ref{sec:motivating-example}, for which we demonstrated that neither the strict projection lemma nor Lemma~\ref{lem:ns-proj-Helmersson} could be applied. Using Theorem~\ref{thm:ns-proj-lem} we derive a non-strict version of~\eqref{eq:example-strict-cond} that can be used to test marginal stability of the system~\eqref{eq:system}. The resulting condition is presented in~\ref{item:prop-stab-conds-3}, below, along with its counterpart that can be used for observer synthesis in~\ref{item:prop-stab-conds-2}.
\begin{proposition}\label{prop:stab-conds}
    The following statements are equivalent:
    \begin{enumerate}[labelindent=0pt,labelwidth=\widthof{\ref{item:prop-stab-conds-1}},label=(P\ref{prop:stab-conds}.\arabic*),itemindent=0em,leftmargin=!]
        \item \label{item:prop-stab-conds-1} The system~\eqref{eq:system} is marginally stable.
        \item \label{item:prop-stab-conds-2} There exist a symmetric positive-definite matrix $P\in\mathbb{S}^{n_x}_{\succ 0}$ and a matrix $X\in\mathbb{R}^{n_x\times n_x}$ such that
        \begin{equation}
            \begin{bmatrix}
                P & A^\top X^\top\\
                \star & X+X^\top - P
            \end{bmatrix} \succcurlyeq 0.
            \label{eq:stab-conds-2}
        \end{equation}
        \item \label{item:prop-stab-conds-3} There exist a symmetric positive-definite matrix $S\in\mathbb{S}^{n_x}_{\succ 0}$ and a matrix $X\in\mathbb{R}^{n_x\times n_x}$ such that
        \begin{equation}
            \begin{bmatrix}
                S & AX\\
                \star & X+X^\top-S
            \end{bmatrix}\succcurlyeq 0.
            \label{eq:stab-conds-3}
        \end{equation}
    \end{enumerate}
    Moreover, any matrix $X$ satisfying~\eqref{eq:stab-conds-2} or~\eqref{eq:stab-conds-3} is non-singular.
\end{proposition}
\noindent Although this is not the main focus of this note, let us elaborate somewhat on the two characterizations of marginal stability obtained in Proposition~\ref{prop:stab-conds}. They are equivalent, but each is useful in different applications.
For instance,~\ref{item:prop-stab-conds-2} is useful for observer synthesis by replacing $A$ with $A+LC$ and applying the linearizing change of variables $Z_L=XL$. Similarly,~\ref{item:prop-stab-conds-3} can be used for controller synthesis by substituting $A$ with $A+BK$ and performing the linearizing change of variables $KX = Z_K$. Since no products between $A$ and $P$ or $S$ appear, both~\ref{item:prop-stab-conds-2} and~\ref{item:prop-stab-conds-3} can be used, following the same development as in~\cite{Oliveira1999}, to synthesize polytopic/switched observers/controllers for polytopic or switched linear systems, while guaranteeing marginal stability of the closed-loop (or estimation error) system using a polytopic/switched weak Lyapunov function.

\subsection{Interpolation and the matrix S-procedure}
We can also apply the NSPL to derive the interpolation result below, which is a version of~\cite[Lemma A.2]{Scherer2005} with slightly weaker assumptions.
\begin{lemma}\label{lem:interpolation}
    Let $R\in\mathbb{S}_{\prec 0}^{m}$ and
    \begin{equation*}
        P=\begin{bmatrix}
        Q & S\\
        \star & R
    \end{bmatrix}\in\mathbb{S}^{n+m}
    \end{equation*}
    with $Q-SR^{-1}S^\top \succcurlyeq 0$. Then, for any $z\in\mathbb{R}^{n}$ and $w\in\mathbb{R}^{m}$, there exists some matrix $\Delta\in\mathbb{R}^{m\times n}$ such that
    \begin{equation}\label{eq:interpolation-Delta}
        w = \Delta z\quad\text{and}\quad\begin{bmatrix}
        I\\
        \Delta
        \end{bmatrix}^\top P\begin{bmatrix} I\\
        \Delta\end{bmatrix}\succcurlyeq 0,
    \end{equation}
    if and only if
    \begin{equation}
        \begin{bmatrix}
            z\\
            w
        \end{bmatrix}^\top P\begin{bmatrix} z\\ w\end{bmatrix} \geqslant 0.\label{eq:interpolation-P}
    \end{equation}
\end{lemma}
\noindent Interpolation results, such as the one in Lemma~\ref{lem:interpolation}, have been used in LMI relaxations techniques for robust control~\cite{Scherer2005}. Interestingly, Lemma~\ref{lem:interpolation} can also be used to derive the matrix S-lemma below.
\begin{lemma}\label{lem:matrix-S-lemma}
    Let $M\in\mathbb{S}^{n+m}$ and
    \begin{equation*}
        N=\begin{bmatrix} N_{11} & N_{12}\\ \star & N_{22}\end{bmatrix}\in\mathbb{S}^{n+m},
    \end{equation*}
    with $N_{22}\prec 0$ and $N_{11}-N_{12}N_{22}^{-1}N_{12}^\top\succcurlyeq 0$. Then, the following statements are equivalent:
    \begin{enumerate}[labelindent=0pt,labelwidth=\widthof{\ref{item:matrix-S-lemma-2}},label=(L\ref{lem:matrix-S-lemma}.\arabic*),itemindent=0em,leftmargin=!]
        \item \label{item:matrix-S-lemma-1} $\begin{bmatrix}
        I\\
        Z
        \end{bmatrix}^{\top} M\begin{bmatrix}I\\ Z\end{bmatrix}\succ 0$ for all $Z$ with $\begin{bmatrix} I\\ Z\end{bmatrix}^{\top} N\begin{bmatrix} I\\ Z\end{bmatrix}\succcurlyeq 0$.
        \item \label{item:matrix-S-lemma-2} There exists some $\alpha\geqslant 0$ such that $M-\alpha N\succ 0$.
    \end{enumerate}
\end{lemma}
\noindent Lemma~\ref{lem:matrix-S-lemma} is~\cite[Corollary 12]{vanWaarde2022} without the requirement that $N$ is non-singular. Similar non-strict results have proved instrumental in recent (direct) data-driven control applications, see, e.g.,~\cite{vanWaarde2022,vanWaarde2022b,Berberich2022,Bisoffi2021,Bisoffi2022}, which, in turn, demonstrates the relevance of the NSPL also in this area.

\subsection{Matrix dilations}
Finally, we indicate that Theorem~\ref{thm:ns-proj-lem} also has applications in matrix (or operator) dilation theory~\cite{Davis1982}, such as in the result stated below.
\begin{lemma}\label{lem:dilation}
    Let $A\in\mathbb{R}^{m\times n}$, $B\in\mathbb{R}^{m\times p}$ and $C\in\mathbb{R}^{q\times n}$. Then, there exists $D\in\mathbb{R}^{q\times p}$ with
    \begin{equation}
        \left\|\begin{bmatrix}
            A & B\\
            C & D
        \end{bmatrix}\right\|\leqslant 1
        \label{eq:dilated}
    \end{equation}
    if and only if
    \begin{equation}
        \left\|\begin{bmatrix} A & B\end{bmatrix}\right\|\leqslant 1\text{ and }\left\|\begin{bmatrix}
            A\\
            C
        \end{bmatrix}\right\|\leqslant 1.\label{eq:dilation-conds}
    \end{equation}
\end{lemma}
\noindent A proof based on Theorem~\ref{thm:ns-proj-lem} can be found in the Appendix.

\section{Conclusions}\label{sec:conclusions}
In this technical note, we presented a non-strict generalization of the projection lemma. This non-strict projection lemma was shown to include both the strict projection lemma and an earlier non-strict version of the projection lemma as special cases, thereby showing that our contribution generalizes these existing results. In addition, we showed several applications of this novel non-strict projection lemma, for which existing results could not be applied. One such application is analyzing marginal stability (or performing marginal stabilization) of discrete-time LTI systems, for which we derived several LMI-based conditions. The resulting stability conditions are such that they can be used for controller and observer synthesis. They may be further extended to accommodate synthesis for polytopic/switched linear systems using a polytopic/switched weak Lyapunov function to guarantee marginal stability for the corresponding closed-loop system. We also demonstrate applications of our results to dilation theory as well as interpolation, where, for the latter, we show that the matrix S-lemma, which proves instrumental in the context of (direct) data-driven control, naturally follows.

\bibliographystyle{IEEEtran}
\bibliography{phd-bibtex}

\begin{thebibliography}{10}
\providecommand{\url}[1]{#1}
\csname url@samestyle\endcsname
\providecommand{\newblock}{\relax}
\providecommand{\bibinfo}[2]{#2}
\providecommand{\BIBentrySTDinterwordspacing}{\spaceskip=0pt\relax}
\providecommand{\BIBentryALTinterwordstretchfactor}{4}
\providecommand{\BIBentryALTinterwordspacing}{\spaceskip=\fontdimen2\font plus
\BIBentryALTinterwordstretchfactor\fontdimen3\font minus
  \fontdimen4\font\relax}
\providecommand{\BIBforeignlanguage}[2]{{%
\expandafter\ifx\csname l@#1\endcsname\relax
\typeout{** WARNING: IEEEtran.bst: No hyphenation pattern has been}%
\typeout{** loaded for the language `#1'. Using the pattern for}%
\typeout{** the default language instead.}%
\else
\language=\csname l@#1\endcsname
\fi
#2}}
\providecommand{\BIBdecl}{\relax}
\BIBdecl

\bibitem{Boyd1994}
S.~{Boyd}, L.~{El Ghaoui}, E.~{Feron}, and V.~{Balakrishnan}, \emph{Linear
  matrix inequalities in control theory}, ser. Studies in Applied
  Mathematics.\hskip 1em plus 0.5em minus 0.4em\relax SIAM, 1994, vol.~15.

\bibitem{Gahinet1994}
P.~M. {Gahinet} and P.~{Apkarian}, ``A linear matrix inequality approach to
  {$H_{\infty}$} control,'' \emph{Int. J. Robust Nonlinear Control}, vol.~4,
  no.~4, pp. 421--448, 1994.

\bibitem{Scherer2000}
C.~W. {Scherer} and S.~{Weiland}, ``Linear matrix inequalities in control,''
  2000, {Lecture} Notes, Dutch Institute for Systems and Control, Delft, the
  Netherlands.

\bibitem{Feng2013}
Y.~{Feng} and M.~{Yagoubi}, ``On state feedback {$H_{\infty}$} control for
  discrete-time singular systems,'' \emph{IEEE Trans. Autom. Control}, vol.~58,
  no.~10, pp. 2674--2679, 2013.

\bibitem{Packard1991}
A.~{Packard}, K.~{Zhou}, P.~{Pandey}, and G.~{Becker}, ``A collection of robust
  control problems leading to {LMI's},'' in \emph{30th IEEE Conf. Decis.
  Control}, 1991, pp. 1245--1250.

\bibitem{Skelton1998}
R.~E. {Skelton}, T.~{Iwasaki}, and K.~M. {Grigoriadis}, \emph{A unified
  algebraic approach to linear control design}.\hskip 1em plus 0.5em minus
  0.4em\relax Taylor \& Francis, 1998.

\bibitem{Packard1994}
A.~{Packard}, ``Gain scheduling via linear fractional transformations,''
  \emph{Syst. Control Lett.}, vol.~22, no.~2, pp. 79--92, 1994.

\bibitem{Apkarian1995}
P.~{Apkarian} and P.~{Gahinet}, ``A convex characterization of gain-scheduled
  {$\mathcal{H}_{\infty}$} controllers,'' \emph{IEEE Trans. Autom. Control},
  vol.~40, no.~5, pp. 853--864, 1995.

\bibitem{Oliveira2002}
M.~C. {de Oliveira}, J.~C. {Geromel}, and J.~{Bernussou}, ``Extended {$H_2$}
  and {$H_{\infty}$} norm characterizations and controller parametrizations for
  discrete-time systems,'' \emph{Int. J. Control}, vol.~75, pp. 666--679, 2002.

\bibitem{Ebihara2015}
Y.~{Ebihara}, D.~{Peaucelle}, and D.~{Arzelier}, \emph{{S}-variable approach to
  {LMI}-based robust control}.\hskip 1em plus 0.5em minus 0.4em\relax Springer,
  2015.

\bibitem{Helmersson1995}
A.~Helmersson, ``Methods for gain scheduling,'' Ph.D. dissertation,
  Link\"{o}ping University, 1995.

\bibitem{Scherer2005}
C.~W. {Scherer}, ``Relaxations for robust linear matrix inequality problems
  with verifications for exactness,'' \emph{SIAM J. Matrix Anal. Appl.},
  vol.~27, no.~2, pp. 365--395, 2005.

\bibitem{vanWaarde2022}
H.~J. {van Waarde}, M.~K. {Camlibel}, and M.~{Mesbahi}, ``From noisy data to
  feedback controllers: {Nonconservative} design via a matrix {S}-lemma,''
  \emph{IEEE Trans. Autom. Control}, vol.~67, no.~1, pp. 162--175, 2022.

\bibitem{Saberi2012}
A.~{Saberi}, A.~A. {Stoorvogel}, and P.~{Sannuti}, \emph{Internal and external
  stabilization of linear systems with constraints}.\hskip 1em plus 0.5em minus
  0.4em\relax Birkh\"{a}user, 2012.

\bibitem{Kellett2004}
C.~M. {Kellett} and A.~R. {Teel}, ``Weak converse {Lyapunov} theorems and
  control-{Lyapunov} functions,'' \emph{SIAM J. Control Optim.}, vol.~42,
  no.~6, pp. 1934--1959, 2004.

\bibitem{Khalil1996}
H.~K. {Khalil}, \emph{Nonlinear Systems}, 3rd~ed.\hskip 1em plus 0.5em minus
  0.4em\relax Prentice Hall, 1996.

\bibitem{Oliveira1999}
M.~C. {de Oliveira}, J.~{Bernussou}, and J.~C. {Geromel}, ``A new discrete-time
  robust stability condition,'' \emph{Syst. Control Lett.}, vol.~37, no.~4, pp.
  261--265, 1999.

\bibitem{Heemels2010}
W.~P.~M.~H. {Heemels}, J.~{Daafouz}, and G.~{Millerioux}, ``Observer-based
  control of discrete-time {LPV} systems with uncertain parameters,''
  \emph{IEEE Trans. Autom. Control}, vol.~55, no.~9, pp. 2130--2135, 2010.

\bibitem{Daafouz2001}
J.~{Daafouz} and J.~{Bernussou}, ``Parameter dependent {Lyapunov} functions for
  discrete time systems with time varying parametric uncertainties,''
  \emph{Syst. Control Lett.}, vol.~43, pp. 355--359, 2001.

\bibitem{Oliveira2001}
M.~C. {de Oliveira} and R.~E. {Skelton}, ``Stability tests for constrained
  linear systems,'' in \emph{Perspectives in Robust Control}, 2001, pp.
  241--257.

\bibitem{Apkarian2001}
P.~{Apkarian}, H.~D. {Tuan}, and J.~{Bernussou}, ``Continuous-time analysis,
  eigenstructure assignment, and {$H_2$} synthesis with enhanced {LMI}
  conditions,'' \emph{IEEE Trans. Autom. Control}, vol.~46, no.~12, pp.
  1941--1946, 2001.

\bibitem{vanWaarde2022b}
H.~J. {van Waarde}, M.~K. {Camlibel}, P.~{Rapisarda}, and H.~L. {Trentelman},
  ``Data-driven dissipativity analysis: {Application} of the matrix
  {S}-lemma,'' \emph{IEEE Control Syst.}, vol.~42, no.~3, pp. 140--149, 2022.

\bibitem{Berberich2022}
J.~{Berberich}, C.~W. {Scherer}, and F.~{Allg\"{o}wer}, ``Combining prior
  knowledge and data for robust controller design,'' \emph{ÏEEE Trans. Autom.
  Control}, pp. 1--16, 2022.

\bibitem{Bisoffi2021}
A.~{Bisoffi}, C.~{De Persis}, and P.~{Tesi}, ``Trade-offs in learning
  controllers from noisy data,'' \emph{Syst. Control Lett.}, vol. 154, p.
  104985, 2021.

\bibitem{Bisoffi2022}
------, ``Data-driven control via {Peterson's} lemma,'' \emph{Automatica}, vol.
  145, p. 110537, 2022.

\bibitem{Davis1982}
C.~{Davis}, W.~M. {Kahan}, and H.~F. {Weinberger}, ``Norm-preserving dilations
  and their applications to optimal error bounds,'' \emph{SIAM J. Numer.
  Anal.}, vol.~19, no.~3, pp. 445--469, 1982.

\bibitem{Zi-Zong2010}
Y.~{Zi-Zong} and G.~{Jin-Hai}, ``Some equivalent results with {Yakubovich}'s
  {S}-lemma,'' \emph{SIAM J. Control Optim.}, vol.~48, no.~7, pp. 4474--4480,
  2010.

\bibitem{Finsler1937}
P.~{Finsler}, ``\"{Uber} das vorkommen definiter und semidefiniter formen in
  scharen quadratischer formen,'' \emph{Comment. Math. Helv.}, vol.~9, pp.
  188--192, 1936--1937.

\end{thebibliography}

\section*{Appendix}
\subsection{Preliminaries}
\begin{lemma}[{Schur complement~\cite[p.~8, 28]{Boyd1994}}]\label{lem:ns-schur}
    Let $Q\in\mathbb{H}^{m}$, $R\in\mathbb{H}^{n}$ and $S\in\mathbb{C}^{m\times n}$. Then,
    \begin{equation}
        \begin{bmatrix}
            Q & S\\
            \star & R
        \end{bmatrix}\succcurlyeq 0\label{eq:schur-mat}
    \end{equation}
    if and only if $R\succcurlyeq 0$, $Q-SR^+ S^{\hop }\succcurlyeq 0$ and $S(I-RR^+)=0$. If $R$ is non-singular,~\eqref{eq:schur-mat} holds if and only if $R\succ 0$ and $Q-SR^{-1}S^{\hop}\succcurlyeq 0$.
\end{lemma}

\begin{lemma}[{S-lemma~\cite[p.~24]{Boyd1994}}]\label{lem:S-lemma}
    Let $M,N\in\mathbb{S}^{n}$ and suppose that there exists some $\bar{x}\in\mathbb{R}^{n}$ such that $\bar{x}^\top N\bar{x}>0$. Then, the following statements are equivalent:
    \begin{enumerate}[labelindent=0pt,labelwidth=\widthof{\ref{item:S-lemma-2}},label=(L\ref{lem:S-lemma}.\arabic*),itemindent=0em,leftmargin=!]
        \item \label{item:S-lemma-1} $x^\top Mx>0$ for all $x\in\mathbb{R}^{n}\setminus\{0\}$ such that $x^\top Nx\geqslant 0$.
        \item \label{item:S-lemma-2} There exists $\alpha\geqslant 0$ such that $M-\alpha N\succ 0$.
    \end{enumerate}
\end{lemma}

\begin{lemma}[{Finsler's lemma~\cite[Theorem 2.2]{Zi-Zong2010},~\cite{Finsler1937}}]\label{lem:finsler}
	Let $M,N\in\mathbb{S}^{n}$. Then, the following statements are equivalent:
    \begin{enumerate}[labelindent=0pt,labelwidth=\widthof{\ref{item:finsler-2}},label=(L\ref{lem:finsler}.\arabic*),itemindent=0em,leftmargin=!]
        \item \label{item:finsler-1} $x^\top Mx>0$ for all $x\in\mathbb{R}^{n}\setminus\{0\}$ such that $x^\top Nx=0$.
        \item \label{item:finsler-2} There exists $\alpha\in\mathbb{R}$ such that $M-\alpha N\succ 0$.
    \end{enumerate}
\end{lemma}

\subsection{Proof of Theorem~\ref{thm:ns-proj-lem}}
    {\bf Necessity: } Suppose there exists $X\in\mathbb{C}^{m\times n}$ such that~\eqref{eq:ns-proj-lem-X} holds. By definition of the annihilators $U_\perp$ and $V_\perp$, we have $UU_\perp=0$ and $VV_\perp=0$. Hence, we conclude
    \begin{equation*}
        U_{\perp}^{\hop } \left(Q+U^{\hop}XV+V^{\hop}X^{\hop}U\right)U_{\perp}=U_{\perp}^{\hop}QU_{\perp}\succcurlyeq 0
    \end{equation*}
    and
    \begin{equation*}
        V_{\perp}^{\hop} \left(Q+U^{\hop}XV+V^{\hop}X^{\hop}U\right)V_{\perp}=V_{\perp}^{\hop}QV_{\perp}\succcurlyeq 0,
    \end{equation*}
    which shows that~\eqref{eq:ns-proj-lem-kernels} holds. It remains to show that~\eqref{eq:ns-proj-coupling} holds as well. Let $x\in\ker U\cap\ker V\cap\{\xi\in\mathbb{C}^{p}\mid \xi^\hop Q\xi=0\}$. If  $S\in\mathbb{H}_{\succcurlyeq 0}^{p}$ is a matrix with $S^2=Q+U^\hop XV + V^\hop X^\hop U$, then
    \begin{equation*}
        \|Sx\|^2 = x^\hop (Q+U^\hop XV + V^\hop X^\hop U)x=x^\hop Qx=0.
    \end{equation*}
    We obtain
    \begin{equation*}
        0 = Sx = S^2x = (Q+U^\hop XV + V^\hop X^\hop U)x=Qx,
    \end{equation*}
    i.e., $x\in\ker Q$. This shows that~\eqref{eq:ns-proj-coupling} holds.

    {\bf Sufficiency:} Suppose that~\eqref{eq:ns-proj-lem-kernels} and~\eqref{eq:ns-proj-coupling} hold. Let $T\in\mathbb{C}^{p\times p}$ be a non-singular matrix, whose columns in the partition $T=\begin{bmatrix} T_1 & T_2 & T_3 & T_4 & T_5\end{bmatrix}$ are chosen to satisfy
   \begin{align}
        \im T_4 &= \ker U\cap \ker V\cap \ker Q,\label{eq:imS4}\\
        \im\begin{bmatrix} T_3 & T_4\end{bmatrix} &= \ker U\cap\ker V,\\
        \im \begin{bmatrix} T_1 & T_3 & T_4\end{bmatrix} &= \ker U,\label{eq:imS1S3S4}\\
        \im\begin{bmatrix} T_2 & T_3 & T_4\end{bmatrix} &= \ker V.\label{eq:imS2S3S4}
    \end{align}
    By a congruence transformation with $T$,~\eqref{eq:ns-proj-lem-X}
    translates into
    \begin{equation}
        Y\coloneqq T^\hop QT + (U T)^\hop X(VT) + (VT)^\hop X^\hop (UT) \succcurlyeq 0.
        \label{eq:Y}
    \end{equation}
    Here, the lack of strict inequalities necessitated a more involved congruence transformation (with $T$) if compared to the SPL's proof. Furthermore, the SPL is proved by constructing some $X$ such that 
    $$\begin{bmatrix} T_1 & T_2 & T_3\end{bmatrix}^\hop (Q+U^\hop XV+V^\hop X^\hop U)\begin{bmatrix} T_1 & T_2 & T_3\end{bmatrix}\succ 0.$$ 
    In contrast, we will see that such a strict inequality cannot be achieved in the current proof. As a result, we have to carefully determine $X$ such that $$\begin{bmatrix}T_1& T_2 & T_3\end{bmatrix}^\hop (Q+U^\hop XV+V^\hop X^\hop U)T_5,$$ satisfies suitable nullspace properties.

    Using~\eqref{eq:imS4}, let us proceed by partitioning $W\coloneqq T^\hop QT$ according to $T$ as
    \begin{equation}
        W = \left[\begin{array}{@{}ccc;{2pt/2pt}c;{2pt/2pt}c@{}}
            W_{11} & W_{12} & W_{13} & 0 & W_{15}\\
            \star & W_{22} & W_{23} & 0 & W_{25}\\
            \star & \star & W_{33} & 0 & W_{35}\\\hdashline[2pt/2pt]
            0 & 0 & 0 & 0 & 0\\\hdashline[2pt/2pt]
            \star & \star & \star & 0 & W_{55}
        \end{array}\right].
        \label{eq:W}
    \end{equation}
    Similarly, using~\eqref{eq:imS1S3S4},~\eqref{eq:imS2S3S4} and~\eqref{eq:imS4}, the term $(UT)^\hop X(VT)$ in~\eqref{eq:Y} reads as
    \begin{equation}
        (UT)^\hop X(VT) = \begin{bmatrix}
            0\\
            (UT_2)^\hop\\
            0\\
            0\\
            (UT_5)^\hop
        \end{bmatrix}X\begin{bmatrix} (VT_1)^\hop\\ 0\\ 0\\ 0\\ (VT_5)^\hop\end{bmatrix}^\hop.
        \label{eq:TGammaXOmegaT}
    \end{equation}
    It follows from~\eqref{eq:imS1S3S4} and~\eqref{eq:imS2S3S4}, respectively, that $\begin{bmatrix}UT_2 & UT_5\end{bmatrix}$ and $\begin{bmatrix} VT_1 & VT_5\end{bmatrix}$ have full column rank. Using~\eqref{eq:W} and~\eqref{eq:TGammaXOmegaT}, $Y$ in~\eqref{eq:Y} reads as
    \begin{align}
        &\left[\begin{array}{@{}c;{2pt/2pt}c;{2pt/2pt}c@{}}
            Y_1 & 0 & Y_2\\\hdashline[2pt/2pt]
            0 & 0 & 0\\\hdashline[2pt/2pt]
            \star & 0 & Y_3
        \end{array}\right] = \label{eq:cond-KLMN}\\
        &\left[\begin{array}{@{}ccc;{2pt/2pt}c;{2pt/2pt}c@{}}
            W_{11} & W_{12}+K^\hop & W_{13} & 0 & W_{15}+M^\hop\\
            \star & W_{22} & W_{23} & 0 & W_{25}+L\\
            \star & \star & W_{33} & 0 & W_{35}\\\hdashline[2pt/2pt]
            0 & 0 & 0 & 0 & 0\\\hdashline[2pt/2pt]
            \star & \star & \star & 0 & W_{55}+ N+N^\hop
        \end{array}\right]\succcurlyeq 0,\nonumber
    \end{align}
    where
    \begin{equation*}
        \begin{bmatrix}
            K & L\\
            M & N
        \end{bmatrix} = \begin{bmatrix}
            (UT_2)^\hop\\
            (UT_5)^\hop
        \end{bmatrix}X\begin{bmatrix}
            VT_1 & VT_5
        \end{bmatrix}.
    \end{equation*}
    Since $\begin{bmatrix} UT_2 & UT_5\end{bmatrix}$ and $\begin{bmatrix} VT_1 & VT_5\end{bmatrix}$ have full column rank, observe that
    \begin{equation}
        X = \begin{bmatrix}
            (UT_2)^\hop\\
            (UT_5)^\hop
        \end{bmatrix}^+\begin{bmatrix}
            K & L\\
            M & N
        \end{bmatrix}\begin{bmatrix}
            VT_1 & VT_5
        \end{bmatrix}^+
        \label{eq:X-reconstr}
    \end{equation}
    satisfies~\eqref{eq:ns-proj-lem-X} for any $K$, $L$, $M$ and $N$ that satisfy~\eqref{eq:cond-KLMN}. In the remainder of this proof, we construct such $K$, $L$, $M$ and $N$.

    First, we construct $K$ that renders $Y_1$ in~\eqref{eq:cond-KLMN} positive semi-definite. To this end, note that, due to~\eqref{eq:ns-proj-lem-kernels},~\eqref{eq:imS1S3S4} and~\eqref{eq:imS2S3S4},
    \begin{equation}
        \begin{bmatrix}
            W_{11} & W_{13}\\
            \star & W_{33}
        \end{bmatrix}\succcurlyeq 0\text{ and }\begin{bmatrix}
            W_{22} & W_{23}\\
            \star & W_{33}
        \end{bmatrix}\succcurlyeq 0.\label{eq:ker-conds-W}
    \end{equation}
    It also follows from~\eqref{eq:ns-proj-coupling} that
    \begin{equation*}
        x^\hop Qx \neq 0\text{ for all }x\in\ker U\cap \ker V\text{ with }x\notin\ker Q.
    \end{equation*}
    Thus, by construction of $T_3$ and using~\eqref{eq:ker-conds-W}, we infer
        \begin{equation*}
        W_{33} = T_3^\hop QT_3 \succ 0.
    \end{equation*}
    Hence, we can apply Lemma~\ref{lem:ns-schur} to~\eqref{eq:ker-conds-W} to obtain
    \begin{equation}
        W_{11} - W_{13}W_{33}^{-1}W_{13}^\hop \succcurlyeq 0\text{ and }W_{22}-W_{23}W_{33}^{-1}W_{23}^\hop \succcurlyeq 0.
        \label{eq:ker-conds-schured}
    \end{equation}
    Since $W_{33}\succ 0$, Lemma~\ref{lem:ns-schur} reveals that $Y_1\succcurlyeq 0$ if and only if
    \begin{equation*}
        \begin{bmatrix}
            W_{11} & W_{12}+K^\hop\\
            \star & W_{22}
        \end{bmatrix} - \begin{bmatrix} W_{13}\\ W_{23} \end{bmatrix} W_{33}^{-1}\begin{bmatrix}W_{13}\\ W_{23} \end{bmatrix}^\hop \succcurlyeq 0.
    \end{equation*}
    By~\eqref{eq:ker-conds-schured}, $K = -W_{12}^\hop + W_{23}W_{33}^{-1}W_{13}^\hop$ renders the latter inequality valid and, hence, $Y_1\succcurlyeq 0$.

    Next, we apply Lemma~\ref{lem:ns-schur} to see that~\eqref{eq:cond-KLMN} is equivalent to
    \begin{equation*}
        Y_1\succcurlyeq 0,~Y_3-Y_2^\hop Y_1^+ Y_2\succcurlyeq 0\text{ and }Y_2^\hop (I-Y_1Y_1^+)=0.
    \end{equation*}
    We have already constructed a matrix $K$ such that $Y_1\succcurlyeq 0$. Let us now construct $L$ and $M$ such that $Y_2^\hop(I-Y_1Y_1^+)=0$, which, due to the symmetry of $Y_1Y_1^+$, is equivalent to $(I-Y_1Y_1^+)Y_2$. Hence, it suffices to find $L$ and $M$ such that we can write $Y_2$ as $Y_2 = Y_1\tilde{P}$ for some $\tilde{P}$, i.e.,
    \begin{equation}
        \begin{bmatrix}
            W_{15} + M^\hop\\
            W_{25} + L\\
            W_{35}
        \end{bmatrix} = \begin{bmatrix}
            W_{11} & W_{13}W_{33}^{-1}W_{23}^\hop & W_{13}\\
            \star & W_{22} & W_{23}\\
            \star & \star & W_{33}
        \end{bmatrix}\tilde{P},
        \label{eq:LMPcond}
    \end{equation}
    where we have substituted the earlier constructed matrix $K$. A particular choice of $\tilde{P}$, $L$ and $M$ that satisfies~\eqref{eq:LMPcond} is
    \begin{equation*}
        \tilde{P}=\begin{bmatrix}
            0\\
            0\\
            W_{33}^{-1}W_{35}
        \end{bmatrix},~\left[\begin{array}{@{}c@{}}
        M^\hop\\\hdashline[2pt/2pt]
        L
        \end{array}\right]=\left[\begin{array}{@{}c@{}}
            -W_{15} + W_{13} W_{33}^{-1}W_{35}\\\hdashline[2pt/2pt]
            -W_{25}+W_{23}W_{33}^{-1}W_{35}
        \end{array}\right].
    \end{equation*}
    It remains to construct $N$ such that
    \begin{equation}
        0\preccurlyeq Y_3 - Y_2^\hop Y_1^+ Y_2 = W_{55}+N+N^\hop - Y_2^\hop Y_1^+Y_2,\label{eq:tuneN}
    \end{equation}
    which we achieve by choosing $N=\alpha I$ with $\alpha>0$ sufficiently large to ensure that~\eqref{eq:tuneN} holds. Since we have found $K$, $L$, $M$ and $N$ for which~\eqref{eq:cond-KLMN} holds, $X$ as in~\eqref{eq:X-reconstr} satisfies~\eqref{eq:ns-proj-lem-X}, which completes the proof.

\subsection{Proof of Corollary~\ref{cor:SPL}}
	We prove Lemma~\ref{lem:SPL} by relying on Theorem~\ref{thm:ns-proj-lem}.
	
{\bf Necessity:} Suppose that~\eqref{eq:trad-proj-1} holds and let $\tilde{Q}=Q-\epsilon I$ for some $\epsilon>0$ such that $Q+U^\hop XV+V^\hop X^\hop U\succ \epsilon I$. Hence,
\begin{equation*}
    \tilde{Q} +U^\hop XV+V^\hop X^\hop U\succcurlyeq 0.
\end{equation*}
By Theorem~\ref{thm:ns-proj-lem},  we infer $U_{\perp}^\hop \tilde{Q}U_{\perp}\succcurlyeq 0$ and $V_{\perp}^\hop \tilde{Q}V_{\perp}\succcurlyeq 0$. Since $U_\perp$ and $V_\perp$ have full column rank, we conclude
\begin{equation*}
    U_\perp^\hop QU_\perp \succcurlyeq \epsilon U_\perp^\hop U_\perp \succ 0
    \text{ and }
    V_\perp^\hop QV_\perp\succcurlyeq\epsilon V_\perp^\hop V_\perp \succ 0.
\end{equation*}

{\bf Sufficiency:} Suppose that~\eqref{eq:trad-proj-2} holds. We can find some sufficiently small $\epsilon > 0$ such that $U_\perp^\hop QU_\perp \succ \epsilon U_{\perp}^\hop U_\perp$ and $V_\perp^\hop QV_\perp\succ \epsilon V_\perp^\hop V_\perp$. Then, $\tilde{Q}:=Q-\epsilon I$ satisfies
\begin{equation}
    U_\perp^\hop \tilde{Q}U_\perp\succ 0\text{ and }V_\perp^\hop \tilde{Q}V_\perp\succ 0.\label{eq:ker-Q-tilde}
\end{equation}
Thus,~\eqref{eq:ns-proj-lem-kernels} holds. It remains to show~\eqref{eq:ns-proj-coupling}, i.e.,
\begin{equation*}
    \ker U\cap\ker V\cap\{\xi\in\mathbb{C}^{p}\mid \xi^\hop \tilde{Q}\xi=0\}\subset \ker\tilde{Q},
\end{equation*}
such that we can apply Theorem~\ref{thm:ns-proj-lem} and conclude \eqref{eq:trad-proj-1} due to the $\epsilon$ perturbation of $Q$. To show~\eqref{eq:ns-proj-coupling}, let $\xi\in \ker U\cap \ker V$ with $\xi^\hop \tilde{Q}\xi=0$. Then, there exists some $\eta$ with $\xi=U_{\perp}\eta$. It follows that $0=\xi^\hop \tilde{Q}\xi = \eta^\hop U_{\perp}^\hop \tilde{Q}U_{\perp}\eta$, which, due to~\eqref{eq:ker-Q-tilde}, implies $\eta =0$ and thus $\xi=0$. This indeed implies $\xi\in\ker\tilde{Q}$, which concludes the proof of~\eqref{eq:ns-proj-coupling}.

\subsection{Proof of Corollary~\ref{cor:Helmersson}}
	We prove Lemma~\ref{lem:ns-proj-Helmersson} by employing Theorem~\ref{thm:ns-proj-lem}. To this end, it suffices to show that~\eqref{eq:Helmersson-extra-cond} and~\eqref{eq:ns-proj-lem-kernels} imply~\eqref{eq:ns-proj-coupling}. Suppose that~\eqref{eq:Helmersson-extra-cond} and~\eqref{eq:ns-proj-lem-kernels} hold and let $\xi\in\ker U\cap\ker V$ be such that $\xi^\hop Q\xi=0$. Then, there exist $\zeta$ and $\eta$ such that $\xi=U_{\perp}\zeta$ and $\xi=V_{\perp}\eta$. It follows that
    \begin{equation*}
        0=\xi^\hop Q\xi = \zeta^\hop U_{\perp}^\hop QU_{\perp}\zeta = \eta^\hop V_{\perp}^\hop QV_{\perp}\eta.
    \end{equation*}
    Using~\eqref{eq:ns-proj-lem-kernels}, this implies that $U_{\perp}^\hop Q\xi=0$ and $V_{\perp}^\hop Q\xi=0$. Hence,
    $$Q\xi\in \ker U_{\perp}^\hop = (\im U_{\perp})^\perp=(\ker U)^\perp = \im U^\hop$$
    and, similarly, $Q\xi \in\im V^\hop$, i.e.,
    \begin{equation*}
        Q\xi\in\im U^\hop \cap \im V^\hop \stackrel{\eqref{eq:Helmersson-extra-cond}}{=} \{0\}.
    \end{equation*}
    Thereby,~\eqref{eq:ns-proj-coupling} holds.

\subsection{Proof of Proposition~\ref{prop:stab-conds}}
    {\bf\ref{item:prop-stab-conds-2}:} It is well-known, see, e.g.,~\cite[p.~211]{Saberi2012}, that~\ref{item:prop-stab-conds-1} is equivalent to the existence of $P\in\mathbb{S}_{\succ 0}^{n_x}$ satisfying~\eqref{eq:stab-conds-1}, which is equivalent to the existence of $P\in\mathbb{S}^{n_x}_{\succ 0}$ such that~\eqref{eq:ns-proj-lem-kernels} holds with
    \begin{equation*}
        Q = \begin{bmatrix}
            P & 0\\
            0 & -P
        \end{bmatrix},~U_{\perp} = \begin{bmatrix}
            I\\
            0
        \end{bmatrix}\text{ and }V_{\perp} = \begin{bmatrix}
            -I\\
            A
        \end{bmatrix}.
    \end{equation*}
    To complete the key step in this proof, we aim to apply Theorem~\ref{thm:ns-proj-lem}. We can take $U = \begin{bmatrix} 0 & I\end{bmatrix}$ and $V = \begin{bmatrix} A & I\end{bmatrix}$. Next, we will show, for this particular problem, that~\eqref{eq:ns-proj-lem-kernels} implies~\eqref{eq:ns-proj-coupling}. To this end, let $x\in\ker U\cap \ker V$ satisfy $x^\top Qx=0$. Note that $x\in\ker U$ implies that $x = (y^\top,0)^\top$ for some $y\in\mathbb{R}^{n_x}$. Moreover, $(y^\top,0) Q(y^\top,0)^\top = y^\top Py=0$ implies $y=0$, since $P\succ 0$ and, thus, $x = 0 \in \ker Q$. It follows that~\eqref{eq:ns-proj-coupling} is valid as well. By Theorem~\ref{thm:ns-proj-lem}, we find that $P\in\mathbb{S}_{\succ 0}^{n_x}$ satisfies~\eqref{eq:stab-conds-1} if and only if there exists $X\in\mathbb{R}^{n_x\times n_x}$ such that
    \begin{equation}
        \begin{bmatrix}
            P & 0\\
            \star & -P
        \end{bmatrix} + \begin{bmatrix} 0\\
        I\end{bmatrix} X\begin{bmatrix} A & I\end{bmatrix} + \begin{bmatrix} A^\top\\
        I\end{bmatrix} X^\top\begin{bmatrix} 0 & I\end{bmatrix}\succcurlyeq 0.\label{eq:PX}
    \end{equation}
    Thus,~\ref{item:prop-stab-conds-1} and~\ref{item:prop-stab-conds-2} are equivalent. To see that $X$ satisfying~\eqref{eq:stab-conds-2} is non-singular, we note that~\eqref{eq:stab-conds-2} implies that $X+X^\top \succcurlyeq P \succ 0$, which holds only if $X$ is non-singular.

    {\bf\ref{item:prop-stab-conds-3}:} The proof is completed using~\eqref{eq:marg-stab-lyap2} and following the  steps as for~\ref{item:prop-stab-conds-2} with the substitutions $P\leftarrow S$ and $A\leftarrow A^\top$.

\subsection{Proof of Lemma~\ref{lem:interpolation}}
{\bf Necessity:} Suppose there exists $\Delta\in\mathbb{R}^{m\times n}$ for which $w=\Delta z$. It immediately follows that~\eqref{eq:interpolation-P} holds, since
\begin{equation*}
    \begin{bmatrix}
        z\\
        w
    \end{bmatrix}^\top P\begin{bmatrix}
        z\\
        w
    \end{bmatrix} = z^\top \begin{bmatrix}
        I\\
        \Delta
    \end{bmatrix}^\top P\begin{bmatrix}
        I\\
        \Delta
    \end{bmatrix}z \stackrel{\eqref{eq:interpolation-Delta}}{\geqslant} 0.
\end{equation*}

{\bf Sufficiency:} Suppose that~\eqref{eq:interpolation-P} holds. If $z=0$, this inequality reads as $w^\top Rw\geqslant 0$. Due to $R\prec 0$, this implies $w=0$. Then, by choosing $\Delta = -R^{-1}S^\top$, we infer $w=\Delta z$ and, by assumption,
\begin{equation*}
    \begin{bmatrix}
        I\\
        \Delta
    \end{bmatrix}^\top P\begin{bmatrix}
        I\\
        \Delta
    \end{bmatrix} = Q-SR^{-1}S^\top \succcurlyeq 0.
\end{equation*}

Next, suppose $z\neq 0$. Using Lemma~\ref{lem:ns-schur} (Schur complement) and $R\prec 0$, the desired inequality in~\eqref{eq:interpolation-Delta} can be expressed as
\begin{equation}
    \begin{bmatrix}
        Q+S\Delta + \Delta^\top S^\top & \Delta^\top\\
        \star & -R^{-1}
    \end{bmatrix}\succcurlyeq 0.
    \label{eq:interp-schured}
\end{equation}
To guarantee $w=\Delta z$, we must have, for some $H\in\mathbb{R}^{m\times n}$,
\begin{equation}
    \Delta = wz^+ + H(I-zz^+).
    \label{eq:DeltaH}
\end{equation}
It follows, by substituting~\eqref{eq:DeltaH} into~\eqref{eq:interp-schured}, that there exists $\Delta$ satisfying~\eqref{eq:interpolation-Delta}, if and only if there exists $H$ such that
\begin{equation}
    \Psi + U^\top HV + V^\top H^\top U \succcurlyeq 0,\label{eq:HHhop}
\end{equation}
with $U = \begin{bmatrix} S^\top & I\end{bmatrix}$, $V = \begin{bmatrix} I-zz^+ & 0\end{bmatrix}$ and
\begin{equation*}
    \Psi = \begin{bmatrix}
        Q+Swz^+ + (Swz^+)^\top & (wz^+)^\top\\
        \star & -R^{-1}
    \end{bmatrix}.
\end{equation*}
The key step of this proof is done by applying Theorem~\ref{thm:ns-proj-lem} to find that a matrix $H$ satisfying~\eqref{eq:HHhop} exists, if and only if~\eqref{eq:ns-proj-lem-kernels} and~\eqref{eq:ns-proj-coupling} hold (for $Q$ replaced by $\Psi$). We introduce the annihilators
\begin{equation*}
    U_{\perp} = \begin{bmatrix}
        I\\
        -S^\top
    \end{bmatrix}\text{ and }V_{\perp} = \begin{bmatrix}
        z & 0\\
        Rw & I
    \end{bmatrix}.
\end{equation*}
We have $U_{\perp}^\top \Psi U_{\perp}=Q-SR^{-1}S^\top\succcurlyeq 0$, by assumption, and
\begin{align*}
    &V^\top_{\perp}\Psi V_{\perp} =\\
    &\quad\begin{bmatrix}
        z^\top Qz + z^\top Sw + (z^\top Sw)^\top +w^\top Rw & 0\\
        \star & -R^{-1}
    \end{bmatrix} \stackrel{\eqref{eq:interpolation-P}}{\succcurlyeq} 0.\nonumber
\end{align*}
Next, we show that~\eqref{eq:ns-proj-coupling} holds. Since $z\neq 0$, we have
\begin{equation}
    \ker U\cap\ker V = \im\begin{bmatrix}
        Iz\\
        -S^\top z
    \end{bmatrix}=\im\begin{bmatrix} I\\
    -S^\top\end{bmatrix}z.
    \label{eq:interp-ker-intersect}
\end{equation}
It follows from~\eqref{eq:interp-ker-intersect} that there exists $x\in\ker U\cap\ker V\setminus\{0\}$ such that $x^\top \Psi x=0$ if and only if
\begin{equation}
    z^\top \begin{bmatrix}
        I\\
        -S^\top
    \end{bmatrix}^\top \Psi\begin{bmatrix}
        I\\
        -S^\top
    \end{bmatrix}z = \|(Q-SR^{-1}S^\top)^{\frac{1}{2}}z\|^2 = 0,\label{eq:QSRSz}
\end{equation}
where we used the fact that, by assumption, $Q-SR^{-1}S^\top\succcurlyeq 0$. It follows that $x^\top \Psi x=0$ for any $x\in\ker U\cap\ker V$ and, hence,~\eqref{eq:ns-proj-coupling} holds if $\ker U\cap\ker V\subset\ker\Psi$. We have
\begin{equation}
    \Psi\begin{bmatrix}U\\V\end{bmatrix}_{\perp} = \begin{bmatrix}
        Q+Swz^+\\
        wz^++R^{-1}S^\top
    \end{bmatrix}z = \begin{bmatrix}
        Qz+Sw\\
        w+R^{-1}S^\top z
    \end{bmatrix}.\label{eq:PsiGOker}
\end{equation}
From~\eqref{eq:interpolation-P} and the fact that $R\prec 0$, we have
\begin{align*}
    0 \leqslant \begin{bmatrix}
        z\\
        w
    \end{bmatrix}^\top P\begin{bmatrix}
        z\\
        w
    \end{bmatrix} &= z^\top Qz + z^\top Sw + w^\top S^\top z + w^\top Rw,\\
    &= (w+R^{-1}S^\top z)^\top R(w+R^{-1}S^\top z)\leqslant 0,
\end{align*}
and, hence, $w = -R^{-1}S^\top z$. Substitution in~\eqref{eq:PsiGOker} yields
\begin{equation*}
    \Psi\begin{bmatrix}U\\V\end{bmatrix}_{\perp} = \begin{bmatrix}
        (Q-SR^{-1}S^\top)z\\
        0
    \end{bmatrix} \stackrel{\eqref{eq:QSRSz}}{=} 0,
\end{equation*}
such that~\eqref{eq:ns-proj-coupling} holds. We apply Theorem~\ref{thm:ns-proj-lem} to conclude that there exists $H$ satisfying~\eqref{eq:HHhop} and, thus, $\Delta$ in~\eqref{eq:DeltaH} satisfies~\eqref{eq:interpolation-Delta}.

\subsection{Proof of Lemma~\ref{lem:matrix-S-lemma}}
{\bf $\mathbf{\ref{item:matrix-S-lemma-2}\Rightarrow\ref{item:matrix-S-lemma-1}}$:} Suppose that~\ref{item:matrix-S-lemma-2} holds. Then, for all $Z\in\mathbb{R}^{m\times n}$ such that $\begin{bmatrix} I & Z^\top\end{bmatrix}N\begin{bmatrix} I & Z^\top\end{bmatrix}^\top \succcurlyeq 0$, we have
\begin{equation*}
    \begin{bmatrix}
        I\\
        Z
    \end{bmatrix}^\top M\begin{bmatrix}
        I\\
        Z
    \end{bmatrix}\succ \alpha \begin{bmatrix} I\\ Z\end{bmatrix}^\top N\begin{bmatrix} I\\
    Z\end{bmatrix} \succcurlyeq 0.
\end{equation*}

{\bf $\mathbf{\ref{item:matrix-S-lemma-1}\Rightarrow\ref{item:matrix-S-lemma-2}}$:} Suppose that~\ref{item:matrix-S-lemma-1} holds. By Lemma~\ref{lem:interpolation}, for any $x=(z^\top,w^\top)^\top$ such that $x^\top Nx\geqslant 0$, there exists $Z\in\mathbb{R}^{m\times n}$ such that
\begin{equation}
    w = Zz\text{ and }\begin{bmatrix}
        I\\
        Z
    \end{bmatrix}^\top N\begin{bmatrix}
        I\\
        Z
    \end{bmatrix}\succcurlyeq 0.\label{eq:wZz}
\end{equation}
Firstly, we consider the case where $N$ has at least one positive eigenvalue. Take any $x=(z^\top,w^\top)^\top\neq 0$ such that $x^\top Nx\geqslant 0$ and let $Z$ be as in~\eqref{eq:wZz}. It follows, by~\ref{item:matrix-S-lemma-1}, that $\begin{bmatrix} I & Z^\top\end{bmatrix}M\begin{bmatrix} I & Z^\top\end{bmatrix}^\top\succ 0$ and, hence,
\begin{equation}
    0<z^\top\begin{bmatrix} I\\
    Z\end{bmatrix}^\top M\begin{bmatrix}I\\ Z\end{bmatrix}z \stackrel{\eqref{eq:wZz}}{=} x^\top Mx.
    \label{eq:xMx}
\end{equation}
Thus,~\ref{item:S-lemma-1} holds and we can apply Lemma~\ref{lem:S-lemma} to conclude that~\ref{item:matrix-S-lemma-2} holds. Secondly, we consider the case where $N$ has no positive eigenvalues, i.e., $N\preccurlyeq 0$. Take any $x=(z^\top,w^\top)^\top\neq 0$ such that $x^\top Nx=0$ and let $Z$ be as in~\eqref{eq:wZz}. Again, it follows, by~\ref{item:matrix-S-lemma-1}, that $\begin{bmatrix} I & Z^\top\end{bmatrix}M\begin{bmatrix} I & Z^\top\end{bmatrix}^\top\succ 0$. Thus,~\eqref{eq:xMx} holds and we can apply Lemma~\ref{lem:finsler} to conclude that there exists some $\alpha\in\mathbb{R}$ such that $M-\alpha N\succ 0$. If $\alpha\geqslant 0$, we are done. Otherwise, since $N\preccurlyeq 0$, we obtain $\alpha N\succcurlyeq 0$ and, thus,
\begin{equation}
    M - 0\cdot N = M\succ \alpha N\succcurlyeq 0,
\end{equation}
which completes our proof as well.

\subsection{Proof of Lemma~\ref{lem:dilation}}
The norm inequality~\eqref{eq:dilated} can be expressed as an LMI, which, using Lemma~\ref{lem:ns-schur}, can be, equivalently, expressed as
\begin{equation}
    0\preccurlyeq \left[\begin{array}{@{}cc;{2pt/2pt}cc@{}}
        I & 0 & A & B\\
        \star & I & C & D\\\hdashline[2pt/2pt]
        \star & \star & I & 0\\
        \star & \star & \star & I
    \end{array}\right] = Q + U^\top DV + V^\top D^\top U,
    \label{eq:dilated-lmi}
\end{equation}
with
\begin{equation*}
    Q = \left[\begin{array}{@{}cc;{2pt/2pt}cc@{}}
        I & 0 & A & B\\
        \star & I & C & 0\\\hdashline[2pt/2pt]
        \star & \star & I & 0\\
        \star & \star & \star & I
    \end{array}\right],~U = \left[\begin{array}{@{}c@{}}
        0\\ I\\\hdashline[2pt/2pt] 0\\ 0
    \end{array}\right]^\top\text{ and }V = \left[\begin{array}{@{}c@{}}
        0\\ 0\\\hdashline[2pt/2pt] 0\\ I
    \end{array}\right]^\top.
\end{equation*}
We introduce the relevant annihilators
\begin{equation}
    U_{\perp} = \left[\begin{array}{@{}ccc@{}}
        I & 0 & 0\\
        0 & 0 & 0\\\hdashline[2pt/2pt]
        0 & I & 0\\
        0 & 0 & I
    \end{array}\right]\text{ and }V_{\perp} = \left[\begin{array}{@{}ccc@{}}
        I & 0 & 0\\
        0 & I & 0\\\hdashline[2pt/2pt]
        0 & 0 & I\\
        0 & 0 & 0
    \end{array}\right].
    \label{eq:dilation-annihilators}
\end{equation}

{\bf Necessity:} By assumption,~\eqref{eq:dilated-lmi} holds. It follows that
\begin{equation}
    U_{\perp}^\top QU_{\perp} = \begin{bmatrix}\begin{smallmatrix}
        I & A & B\\
        \star & I & 0\\
        \star & \star & I
    \end{smallmatrix}\end{bmatrix} \succcurlyeq 0\text{ and }V_{\perp}^\top QV_{\perp} = \begin{bmatrix}\begin{smallmatrix}
        I & 0 & A\\
        \star & I & C\\
        \star & \star & I
    \end{smallmatrix}\end{bmatrix}\succcurlyeq 0,\label{eq:lmi-norm-ineqs}
\end{equation}
which implies the norm inequalities in~\eqref{eq:dilation-conds}.

{\bf Sufficiency:} Suppose~\eqref{eq:dilation-conds} holds, which is equivalent to~\eqref{eq:lmi-norm-ineqs}. It remains to show that~\eqref{eq:ns-proj-coupling} holds, such that we can apply the NSPL which finishes the proof. To this end, note that~\eqref{eq:lmi-norm-ineqs} implies, using Lemma~\ref{lem:ns-schur}, that
\begin{equation}
    \begin{bmatrix}
        I & A\\
        \star & I
    \end{bmatrix}\succcurlyeq 0,\begin{bmatrix}
        I-BB^\top & A\\
        \star & I
    \end{bmatrix}\succcurlyeq 0\,\text{and}\,\begin{bmatrix}
        I-C^\top C & A\\
        \star & I
    \end{bmatrix}\succcurlyeq 0.
    \label{eq:IABBCC}
\end{equation}
Let $x=(x_1^\top,x_2^\top,x_3^\top,x_4^\top)^\top\in\ker U\cap\ker V$ with $x^\top Qx=0$. Due to~\eqref{eq:dilation-annihilators}, we have $x_2 = 0$ and $x_4 = 0$. It follows that
\begin{equation*}
    0 = x^\top Qx = \begin{bmatrix}
        x_1\\
        x_3
    \end{bmatrix}^\top\begin{bmatrix}
        I & A\\
        \star & I
    \end{bmatrix}\begin{bmatrix}
        x_1\\
        x_3
    \end{bmatrix} \stackrel{\eqref{eq:IABBCC}}{=} \left\|\begin{bmatrix}
        I & A\\
        \star & I
    \end{bmatrix}^{\frac{1}{2}}\begin{bmatrix}
        x_1\\
        x_3
    \end{bmatrix}\right\|^2,
\end{equation*}
which implies that
\begin{equation}
    \begin{bmatrix}
        I & A\\
        \star & I
    \end{bmatrix}\begin{bmatrix}
        x_1\\
        x_3
    \end{bmatrix} = 0.
    \label{eq:0cond}
\end{equation}
Next, we observe that
\begin{equation}
    0\stackrel{\eqref{eq:IABBCC}}{\preccurlyeq}\begin{bmatrix}
        x_1\\
        x_3
    \end{bmatrix}^\top\begin{bmatrix}
        I-BB^\top & A\\
        \star & I
    \end{bmatrix}\begin{bmatrix}
        x_1\\
        x_3
    \end{bmatrix} \stackrel{\eqref{eq:0cond}}{=} -\|B^\top x_1\|^2,
\end{equation}
from which we conclude that $B^\top x_1=0$. Similarly, we infer from~\eqref{eq:IABBCC} that $Cx_3=0$. Combining~\eqref{eq:0cond}, $B^\top x_1=0$ and $Cx_3=0$, we obtain
\begin{equation}
    Qx = \left[\begin{array}{@{}cc;{2pt/2pt}cc@{}}
        I & 0 & A & B\\
        \star & I & C & 0\\\hdashline[2pt/2pt]
        \star & \star & I & 0\\
        \star & \star & \star & I
    \end{array}\right]\left[\begin{array}{@{}c@{}}
        x_1\\
        0\\\hdashline[2pt/2pt]
        x_3\\
        0
    \end{array}\right] = \left[\begin{array}{@{}c@{}}
        x_1 + Ax_3\\
        Cx_3\\\hdashline[2pt/2pt]
        A^\top x_1 + x_3\\
        B^\top x_1
    \end{array}\right] = 0.
\end{equation}
In other words, $x\in\ker Q$ and, thus,~\eqref{eq:ns-proj-coupling} holds. Therefore, we can apply the NSPL to conclude that there exists $D\in\mathbb{R}^{q\times p}$ satisfying~\eqref{eq:dilated}.

\end{document}